\documentclass[11pt]{article}
\usepackage[margin=1in]{geometry}
\usepackage{amsmath}
\usepackage{amsfonts}
\usepackage{amssymb}
\usepackage{mathtools}
\usepackage{mathrsfs}
\usepackage[toc,page]{appendix}
\usepackage[nottoc,numbib]{tocbibind}
\usepackage{hyperref}
\usepackage{pst-node}
\usepackage{tikz-cd} 
\usetikzlibrary{positioning}

\makeatletter
\let\Hy@linktoc\Hy@linktoc@none
\makeatother

\begin{document} 
\begin{center}
\textsc{On the invariant cycle theorem for families of Nori motives}\\
\vspace{10px}
Amir Mostaed \\

\begin{center}
\end{center}
\end{center}\vspace{10px}

\textsc{Abstract.} 
In this paper, we  prove a motivic enhancement of the theorem of the fixed part in Hodge theory due to Deligne. 
In the pure motivic case, this was done for the first time by André in [And96].
Our main result is an extension to the mixed case, which strengthens a result by Arapura [Ara13] and also provides an alternative and simpler proof, in the framework of Nori motives, of a result by Ayoub [Ayo14].

\begin{flushleft}
\end{flushleft}
\vspace{10px}

\begin{flushleft}
\textbf{1 \quad Introduction}

\end{flushleft}

Let $k\hookrightarrow \mathbb{C}$ be a subfield. Recall that the category MM$(k)$ of Nori motives is constructed out of quadruples $(X,Y,n, i)$ where $X$ is a $k$-variety, $Y$ is a closed subvariety of $X$, $i$ is a non-negative integer and $j$ an integer. In [Ara13], Arapura has constructed a relative version (i.e., relative to a fix base scheme $S$ over $k$) which is a tannakian category, we denote MM$(S)$. He also showed that [Ara13, Theorem 6.4.1] the category of relative pure Nori motives is equivalent to the category of André's motives defined in terms of motivated cycles  in [And96]. \\
For $M\in \text{MM}(S)$, one can consider the tannakian subcategory $\langle M \rangle$ generated by  $M$. The Betti realization in this setting is a tensor functor $H_\text{B} : \text{MM}(S) \rightarrow \text{LS}(S^\text{an})$, where $\text{LS}(S^\text{an})$ is
the tannakian category of local systems over the analytic variety $S^\text{an}$.
For any $s\in S$, taking the fiber at $s$ provides a fiber functor, and by tannakian duality, one associates to  $\langle M \rangle$ (resp. $\langle H_\text{B}( M) \rangle$) the algebraic group $G_{\text{mot},s}(M)$ (resp. $G_{\text{mono},s}(M)$). 
This way, $H_\text{B}$ provides a morphism of algebraic groups 
$$ G_{\text{mono},s}(M) \rightarrow G_{\text{mot},s}(M). $$
Indeed, this map is an embedding because the monodromy group $G_{\text{mono},s}(M)$ is the
Zariski closure of the image of the monodromy representation $\pi_1(S^\text{an}, s) \rightarrow
\text{GL}(H_\text{B}(M)_s)$, and thus a subgroup of $G_{\text{mot},s}(M)$.

There is yet another algebraic group onstage, corresponding to the largest constant submotive of $M$, which we denote by $G_{\text{mot},0,s}(M)$.  One has a quotient map $$G_{\text{mot},s}(M) 	\twoheadrightarrow G_{\text{mot},0,s}(M).$$

\begin{flushleft}
Patching this two maps together gives a short exact sequence of algebraic groups:\\
\end{flushleft}

Theorem. \textit{The sequence of algebraic groups} 
\begin{equation} 
0 \rightarrow G_{\text{mono},s}(M) \rightarrow G_{\text{mot},s}(M) \rightarrow G_{\text{mot},0,s}(M) \rightarrow 0 \tag{$\ast$} 
\end{equation}  \textit{is exact.}\\
\pagebreak

This 
theorem is due to Nori, Jossen and Ayoub independently and
in different contexts (only Ayoub, in [Ayo14], has published details, but his motives are
defined differently; nevertheless, by a more recent and quite non-trivial theorem
of Choudhury-Gallauer [CG14], his category is equivalent to Nori’s category). Thus, this theorem is equivalent to Ayoub's result [Ayo14, Proposition 37]. We give a simpler proof in this article.  \\
 To prove the exactness, we proceed in two steps: first we show that the group $G_{\text{mono},s}(M)$
is a normal subgroup of $G_{\text{mot},s}(M)$, and then the induced map $G_{\text{mot},s}(M)/G_{\text{mono},s}(M) \rightarrow G_{\text{mot},0,s}(M)$ is surjective.\\

Acknowledgments. This paper consists of part of my Ph.D. thesis.  I am indebted to Prof. Yves André for his continuous support and guidance. \\

\begin{flushleft}
\textbf{2\quad  Preliminaries}
\end{flushleft}  
We shall use Nori's construction of a $\mathbb{Q}$-linear tannakian category MM$(\text{Spec}\:k)$ of mixed motives over a field $k$ of characteristic zero. Forgetting
the tensor product (which is constructed using Nori’s subtle basic
lemma about cellular decompositions), this abelian category is characterized
by certain universal property with respect to functorial maps from a 
diagram category $\Delta$ consisting of  quadruples $(X,Y,i,j)$ where $X$ is a $k$-variety, $Y$ is a closed subvariety of $X$, $i$ is a non-negative integer and $j$ an integer, 
to abelian categories. The corresponding object in MM$(\text{Spec}\; k)$ is denoted by $h^i(X,Y)(j)$. \\

To classical cohomology theories (Betti, de Rham, étale,...), there correspond
tensor functors, called “realizations”, to categories of finite dimensional vector spaces
over appropriate fields. For instance, de Rham cohomology gives rise to
the de Rham realization $H_\text{dR}: \text{MM}(\text{Spec}\: k) \rightarrow \text{Vec}_\mathbb{Q}$, Betti cohomology (which depends on a fixed  embedding of $k$ in $\mathbb{C}$, assumed to exist) gives rise to Betti realization $H_\text{B}: \text{MM}(\text{Spec}\; k) \rightarrow \text{Vec}_\mathbb{Q}$. There is a canonical comparison isomorphism,
called the “period isomorphism”, between the complexification
of the two realizations: $$ H_\text{dR} \underset{k}\otimes \mathbb{C} \cong  H_\text{B} \underset{\mathbb{Q}}\otimes \mathbb{C}.$$ For brevity, one writes $H_\text{dR}^{i}(X,Y)(j)$ instead of $H_\text{dR}(h^{i}(X,Y)(j))$, and  $H_\text{B}^{i}(X,Y)(j)$ instead of $H_\text{B}(h^{i}(X,Y)(j))$ where we view $i$ as cohomological degree and $j$ as Tate twist. \\

There exists an analogue in the relative context, due to Arapura [Ara13], that is, over a smooth
base $S$, which correponds intuitively to the notion of “family of Nori motives
over $S$”. It is constructed from
pairs $(X, Y )$, but in the relative context, that is, over S. We
write $h^i((X ,Y )/S)(j)$ for the relative analogues of $h^i(X, Y )(j)$, associated to the quadruples $(X\overset{f}\rightarrow S\:, Y, i, j)$ consisting of a quasiprojective family $f:X\rightarrow S$ so that $f$ can be completed to a smooth projective map, a subvariety $Y \subset X$ so that $Y$ together with the boundary is a divisor with normal crossings, and indices $i\in \mathbb{N}$, $j\in \mathbb{Z}$. This is a tannakian category such that its image under a suitable Betti realization lies in the category of locally constant sheaves of $\mathbb{Q}$-vector spaces. We denote this category by MM$(S)$. 

Now, let us start from a Nori motive $M$ on Spec $k(S)$ with $k$ an algebraically
closed subfield of $\mathbb{C}$. One can write $M$ as a
subobject of $h^i(X, Y )(j)$. Then up to replacing $S$
by a dense open subset, one may assume that it extends to a pair $(\textbf{X},\textbf{Y})$
satisfying the same condition but relative to $S$. This process is called
“spreading out”: from a motive $M$ on the generic point, one gets a family
of motives $\textbf{M}$ over $S$, that is, an object of $MM(S)$.
In this relative situation, there is an analog of $H_\text{dR}: \text{MM}(S) \rightarrow \text{
MIC}(S)$, where $\text{MIC}(S)$ is the tannakian category of vector bundles over
$S$ endowed with an integrable connection. Thus for any $M \in \text{MM}(S)$,
$H_\text{dR}(M)$ consists of a vector bundle (the de Rham cohomology of the family
$M$ of Nori motives) together with an integrable connection $\nabla_M$, called
the “Gauss-Manin connection”.\\
Similarly, if $S$ is defined over an algebraically closed subfield $k$ of $\mathbb{C}$,
there is an analogue of Betti realization, $H_\text{B} : \text{MM}(S) \rightarrow \text{LS}(S^\text{an})$, where $\text{LS}(S^\text{an})$ is
the tannakian category of local systems over the analytic variety $S^\text{an}$.\\
 There
is also a comparison isomorphism (a special case of the “Riemann-Hilbert
correspondence”), which goes as follows: if one considers the associated
analytic vector bundle $H_\text{dR}(M)^\text{an}$ with connection $\nabla_
M^\text{an}$ , its local system of
“horizontal sections” Ker $\nabla_
M^\text{an}$ is precisely $H_\text{B}(M)$, and $$H_\text{dR}(M)^\text{an} \cong
H_\text{B}(M) \underset{\mathbb{C}}\otimes \mathcal{O}_{S^\text{an}}$$ as an isomorphism of sheaves on $S^\text{an}$.\\

Now consider the tannakian subcategory $\langle \textbf{M} \rangle$ of MM$(S)$
generated by $\textbf{M}$. Because $\textbf{M}$ is determined by $M$, taking the generic point
$\langle \textbf{M} \rangle \rightarrow \langle M \rangle $ is an equivalence of tannakian categories. After spreading out, one can pick a closed point $s \in S$ and consider the
Betti realization of the fiber at $s$, which gives rise to a fiber functor $\omega_s$, and thus provides
a tannakian group $G_{\text{mot},s}(M)$ such that the category $\langle \textbf{M} \rangle \cong \langle M \rangle $ is equivalent
to Rep$_\mathbb{Q}\: G_{\text{mot},s}(M)$ via $\omega_s$.\\
On the other hand,
one can look at $\langle H_\text{dR}(M) \rangle$, $\langle H_\text{dR}(M)^\text{an} \rangle$, $\langle H_\text{B}(M) \rangle$,
the tannakian subcategories of MIC$(S)$, MIC$(S^\text{an})$, LS$(S^\text{an})$ generated
by $H_\text{dR}(M)$, $H_\text{dR}(M)^\text{an}$, $H_\text{B}(M)$ respectively.\\
Taking the fiber at $s$ as a fiber functor, tannakian duality shows
that these are equivalent to categories of representations of appropriate
linear algebraic groups over $\mathbb{Q}$, which we denote by $G_{\text{diff},s}(M)$, $G_{\text{diff},s}(M^\text{an})$, $G_{\text{mono},s}(M)$ respectively. Because the connection $\nabla_M$ underlying $H_\text{dR}(M)$ is regular, the analytification functor $\langle H_\text{dR}(M) \rangle \rightarrow \langle H_\text{dR}(M)^\text{an} \rangle$
is an equivalence, hence $G_{\text{diff},s}(M)= G_{\text{diff},s}(M^\text{an})$. On the other hand,
by the comparison isomorphism, $\langle H_\text{dR}(M)^\text{an} \rangle \rightarrow \langle H_\text{B}(M) \rangle$ is an equivalence as well, so that
$G_{\text{diff},s}(M^\text{an})=G_{\text{mono},s}(M)$.\\
There is one more tannakian category in this context:
the tannakian subcategory $\langle M \rangle _0$ of $\langle M \rangle $ generated by “constant motives”,
i.e. motives of the form $M \underset{k}\otimes S$. More precisely, motives in the
essential image of the base change functor $$\text{MM}(\text{Spec} \; k) \hookrightarrow \text{MM}(S).$$
Applying tannakian duality, one has the quotient $$G_{\text{mot},s}(M)\twoheadrightarrow G_{\text{mot},0,s}(M).$$

\begin{flushleft}
\textbf{3 \quad The motivic theorem of the fixed part}

\end{flushleft}

Let us discuss the main ingredient in the proof of the theorem.\\

The classical theorem of the fixed part, a.k.a. the global invariant
cycle theorem, is a fundamental result in the topology/homology of families
of algebraic varieties. In its simplest form, it states that for a projective smooth morphism $X \rightarrow S$ over $\mathbb{C}$ and
a base point $s \in S$, $H^i(X_s)^{\pi_1(S,s)}=\text{Im}(\iota^\ast_s: H^i(X) \rightarrow H^i(X_s))$. Here $\iota : X_s \subset X$ is the canonical inclusion of the fibre in the total
space.
Stated like this, this is just a description of the monodromy invariant
subspace of cohomology at the level of vector spaces (no extra
structure involved). But this can be upgraded in the abelian category of Hodge structures: indeed the right hand side is a sub-Hodge structure of $H^i(X_s)$, and hence so is $H^i(X_s)^{\pi_1(S,s)}$. This can be upgraded in a similar way in the abelian category of pure motives, see [And96].\\

We now turn to the theorem of the fixed part in the mixed motivic setting.   
Let us first recall that by definition of MM$(\text{Spec}\; k)$, any object of MM$(\text{Spec}\; k)$ is a subquotient
of an object of the form $h^i(X,Y)(j)$.
There are two very useful (and quite non-trivial) results which improve
on this:
\begin{itemize}
\item[(a)] any object of  MM$(\text{Spec}\; k)$ is a subobject, and not only a subquotient, of
an object of the form $h^i(X,Y)(j)$ [FJ18, Theorem 6.13]. 
\item[(b)] one may assume that $X$ is a complement of a divisor $D$ in a projective smooth $k$-variety $\overline{X}$, and that $Y$ is the restriction to $X$ of a normal crossings divisor $\overline{Y}$ in $\overline{X}$ containing $D$ [Ara13, Theorem 6.1.2].
\end{itemize}

The theorem also works in the context of
pairs $(X,Y)$ as in (b) above, in the relative case (over S): $$H_\text{B}(h^i((X, Y )/S)(j))^{\pi_1(S,s)}_s=H^i((X_s, Y_s)(j))^{\pi_1(S,s)}=\text{Im}(\iota^\ast_s: H^i(X,Y)(j) \rightarrow H^i(X_s,Y_s)(j)).$$

As explained in [Ara13, Lemma 7.2.2], this implies that
$H_\text{B}(h^i((X, Y )/S)(j))^{\pi_1(S,s)}_s$ is the Betti realization of a constant submotive
of $h^i((X, Y )/S)(j)$.\\

In tannakian terms this means the
following. We
start with a family $M\in \text{MM}(S)$ of motives over $S$ and build $\langle M \rangle$. We consider
an object $N \in \langle M \rangle$, the local system $H_\text{B}(N)$ on $S^\text{an}$, and its fibre at $s$,
$H_\text{B}(N)_s$, which is a finite dimensional $\mathbb{Q}$-vector space underlying a representation
of $\pi_1(S,s)$. We are interested in the invariant space $H_\text{B}(N)^{\pi_1(S,s)}_s$. Note that $H_\text{B}(N)_s$ is a representation of $G_{\text{mono},s}(M)$. Thus for any family of Nori motives $N \in \langle M \rangle$ of the form $h^i((X, Y )/S)(j)$, $H_\text{B}(N)^{\pi_1(S,s)}_s$ is stable under $G_{\text{mot},s}(M)$, hence is fiber at $s$ of the Betti realization of a submotive of $N$ over $S$ which is constant, that is, belongs to $\langle M \rangle _0$. \\
Let us now pass from motives of the form $h^i((X, Y )/S)(j)$ to the general case.
It is rather straightforward, using the fact that MM$(S)$ is an abelian category, to treat the case of a submotive of such an object $h^i((X, Y )/S)(j)$ but there is a problem with quotients: invariants in a
quotient are duals of coinvariants of the dual, and one cannot reduce to
invariants except if the monodromy is semisimple. Fortunately, (a refinement of) point (a)
above allows to ignore subquotients and to deal only with subobjects. In conclusion: \\

Theorem (Motivic theorem of the fixed part). For any family of Nori motives $N \in \langle M \rangle$, $H_\text{B}(N)^{\pi_1(S,s)}_s$ is stable under $G_{\text{mot},s}(M)$, hence is fiber at $s$ of the Betti realization of a submotive of $N$ over $S$ which is constant, that is, belongs to $\langle M \rangle _0$.\\

\begin{flushleft}
\textbf{4 \quad Conclusion of the proof of the main theorem}

\end{flushleft}

Let us go back to the sequence ($\ast $). The problem of its exactness
splits into two (successive) subproblems: 
\begin{itemize}
\item[1)] The group $G_{\text{mono},s}(M)$
is a normal subgroup of $G_{\text{mot},s}(M)$.
\item[2)] The induced map $G_{\text{mot},s}(M)/G_{\text{mono},s}(M) \rightarrow G_{\text{mot},0,s}(M)$ is surjective (and hence an isomorphism).
\end{itemize}

For 1), one can use the normality criterion of [And92, lemme 1]: it suffices
to show that for any character $\chi: G_{\text{mono},s}(M) \rightarrow \mathbb{G}_m$, the semi-invariant space $(H_\text{B}(N)^{\pi_1(S,s)}_s)^\chi$ is stable under $G_{\text{mot},s}(M)$. We recall $(H_\text{B}(N)^{\pi_1(S,s)}_s)^\chi$ is the space of elements $v \in H_\text{B}(N)^{\pi_1(S,s)}_s$ such that $h(v)=\chi(h).v$ for every $h \in G_{\text{mono},s}(M)$.\\

Now if all characters of $G_{\text{mono},s}(M)$ are trivial, then via the tannakian
dictionary, the desired properties 1) and 2) become nothing but the motivic theorem of the fixed part. Therefore, exactness of ($\ast$) holds if all
characters of $G_{\text{mono},s}(M)$ are trivial.\\
However, in general, $\mathbb{Q}$-characters of $G_{\text{mono},s}(M)$ may be non-trivial;
“but not so much”: given a character $\chi:G_{\text{mono},s}(M)\rightarrow \mathbb{Q}^\times$ that appears  in a representation of $G_{\text{mot},s}(M)$, the target of $\chi$ is $\mathbb{Z}^\times = \{\pm 1\} $, i.e. $\chi$  is either trivial or of order 2. See [And17b, Lemma C.9]. Thus we assume that $\chi$ is of order exactly 2. By passing to a double  étale covering $S'$ of $S$, we further assume that it is trivial. In this case, the \textit{new} sequence ($\ast$) is exact again by above discussion. But how can
one descend from $S'$ to $S$?\\

\begin{flushleft}
\textbf{5 \quad Weil generic points}

\end{flushleft}

One can  start from a motive $M$ over $k(S)$, where $k$ is an algebraically closed subfield
of $\mathbb{C}$, which spreads out as an object $\mathbf{M}$ of MM$(S)$. Without loss of generality, one may assume that $k$ is countable (since a
motive is defined in terms of pairs of varieties which themselves involve
only finitely many equations, and the same for $S$). This allows to find an
embedding of $k$-extensions $k(S) \xhookrightarrow {} \mathbb{C}$. Such an embedding corresponds to
a complex point $s \in S(\mathbb{C})$ and complex points of this type (i.e. such that
Spec($\mathbb{C}) \rightarrow S$ maps to the generic point of $S$) are called \textit{Weil generic points}.
The point in using Weil generic points is that the realizations at a Weil generic point $s$ make sense both for $\langle M\rangle $ and for $\langle \textbf{M}\rangle$, and moreover they coincide. In particular, $G_{\text{mot},s}(M )\cong G_{\text{mot},s}( \textbf{M})$.
 Therefore, one may  choose a Weil generic point $s$ of $S$ in order to define the Betti realization
$\omega_s$ using $s$ (instead of $H_{B}(\cdot)_s$).   \\ One  can extend the Weil generic point $s$ to an embedding of the algebraic
closure $\overline{k(S)} \xhookrightarrow {} \mathbb{C}$, which thus defines a Weil generic point of every finite
covering of $S$, still denoted by $s$.\\

\begin{flushleft}
\textbf{6 \quad An exact sequence for monodromy groups}

\end{flushleft}

Let $G_{\text{mono},s}(M)$ be again the algebraic monodromy group, i.e. the
Zariski closure of the image of the monodromy representation $\pi_1(S, s) \rightarrow
\text{GL}(\omega_s(\mathbf{M}))$. This is a subgroup of $G_{\text{mot},s}(M)$.\\
By the motivic theorem of the fixed part, the space of $G_{\text{mono},s}(M)$-invariants in $\omega_s(\mathbf{M})$ is the image by $\omega_s$ of a constant subobject of $\mathbf{M}$; in
particular, it is stable under $G_{\text{mot},s}(M)$.\\
The same remains true of course if one replace $M$ by any object of $\langle M \rangle$
and $S$ by any double étale covering $S'$.\\
One has an exact sequence $$ 1 \rightarrow \pi_1(S',s) \rightarrow \pi(S,s) \rightarrow \text{Gal}(k(S')/k(S))=\mathbb{Z}/2 \rightarrow 1,$$ which induces a similar exact sequence for the images of $\pi_1(S,s)$ and $\pi_1(S',s)$ in $\text{GL}(\omega_s(\mathbf{M}))=\text{GL}(\omega_s(\mathbf{M'}))$, and also for their Zariski closures $$ 1 \rightarrow G_{\text{mono},s}(M_{S'}) \rightarrow G_{\text{mono},s}(M) \rightarrow \text{Gal}(k(S')/k(S))=\mathbb{Z}/2  \rightarrow 1.$$ 
\pagebreak

\begin{flushleft}
\textbf{7 \quad An exact sequence for motivic Galois groups}

\end{flushleft}

Our aim is to deduce from the motivic theorem of the fixed part the following
little variation: the space $\omega_s(\mathbf{M})^\chi$ of elements on which $G_{\text{mono},s}(M)$ acts through $\chi$ is stable under $G_{\text{mot},s}(M)$. For this, we may enlarge $G_{\text{mot},s}(M)$, for instance replace it by the motivic
Galois group $G'_{\text{mot},s}(M)$ of the tannakian category generated by $M$ and by the Artin motives defined by the finite extension $k(S')$ of $k(S)$ (a
posteriori, one can see that this is no enlargement). The subcategory of
Artin motives inside this category gives rise to a morphism $G'_{\text{mot},s}(M) \rightarrow \text{Gal}(k(S')/k(S))$, and there is an exact sequence $$ 1 \rightarrow G'_{\text{mot},s}(M_{S'})= G_{\text{mot},s}(M_{S'}) \rightarrow G'_{\text{mot},s}(M) \rightarrow \text{Gal}(k(S')/k(S)) \rightarrow 1,$$ by the work of Jossen [Jos16, Theorem 10.7]. 

\begin{flushleft}
\textbf{8 \quad A commutative diagram of groups}

\end{flushleft}
Now the embedding of monodromy groups in motivic Galois groups is
compatible with the last two exact sequences, so that they
match into a commutative diagram of exact sequences. Since $G_{\text{mot},s}(M_{S'})$ is a normal subgroup of $G'_{\text{mot},s}(M)$, it makes sense to form the compositum
$G_{\text{mot},s}(M_{S'}).G_{\text{mono},s}(M)$. Using the commutative diagram,
this is the full group $G'_{\text{mot},s}(M).$\\

Let $W=\omega_s(\mathbf{M})^\chi \oplus \omega_s(\mathbf{M})^{G_{\text{mono},s}(M)}=\omega_s(\mathbf{M})^{G_{\text{mono},s}(M_{S'})} $. Then $W$ is $G'_{\text{mot},s}(M_{S'})$-stable by the motivic theorem of the fixed part.
It is also stable by $G_{\text{mono},s}(M)$. So $W$ is stable under $G'_{\text{mot},s}(M)$. So we may substitute $W$ to the original representation $\omega_s(\mathbf{M})$; in other
words, we may assume that  $G_{\text{mono},s}(M)$ acts by $\mathbb{Z}/2$ on $\omega_s(\mathbf{M})$, and denoting by
a subscript $\pm$ the $\pm1$-eigenspaces. We have to show that $\omega_s(\mathbf{M})^{-}$ is stable
under $G'_{\text{mot},s}(M).$ Since Hom$_{G_{\text{mono},s}(M)}(\omega_s(\mathbf{M})^+,\omega_s(\mathbf{M})^-)$ = Hom$_{G_{\text{mono},s}(M)}(\omega_s(\mathbf{M})^-,\omega_s(\mathbf{M})^+)=0$, the exact sequence 
\begin{equation}
0\rightarrow \omega_s(\mathbf{M})^+ \rightarrow \omega_s(\mathbf{M}) \rightarrow \omega_s(\mathbf{M})^- \rightarrow 0 \tag{$\star$} 
\end{equation}  has a $\mathbb{Z}/2$-equivariant splitting, and we have to show that this splitting is $G'_{\text{mot},s}(M)$-equivariant. For this, we apply ($\star$) to the dual representation $\omega_s(\mathbf{M})^{\vee}$, and compare this exact sequence to the ($G'_{\text{mot},s}(M)$-equivariant)
dual exact sequence of ($\star$) for $\omega_s(\mathbf{M})$: their $\mathbb{Z}/2$-equivariant splittings
both correspond to the decomposition $\omega_s(\mathbf{M})^{\vee} = (\omega_s(\mathbf{M})^{\vee})^+ \oplus (\omega_s(\mathbf{M})^{\vee})^-$. Therefore the morphisms $$(\omega_s(\mathbf{M})^{\vee})^+ \rightarrow (\omega_s(\mathbf{M})^{+})^{\vee}, \quad (\omega_s(\mathbf{M})/\omega_s(\mathbf{M})^{+})^{\vee} \rightarrow \omega_s(\mathbf{M})^{\vee}/(\omega_s(\mathbf{M})^{\vee})^+ $$ are bijective and we get the claim for $\omega_s(\mathbf{M})^{\vee}$, and hence for $\omega_s(\mathbf{M})$.\\

This proves that the first factor $\omega_s(\mathbf{M})^\chi$ is also stable under $G'_{\text{mot},s}(M)$, and thus under $G_{\text{mot},s}(M)$,
 as desired.\\
 \vspace{20px}

\textit{Remark.} One can consider the category of (relative)\textit{ exponential} motives, along the same construction of category of Nori motives, by introducing \textit{potentials} associated to varieties over complex numbers. Typical objects are of the form $(X\rightarrow S, Y, n, f, i)$ where $X\rightarrow S$ is a morphim of quasi-projective varieties, $Y$ a closed subvariety of $X$, $f\in \mathcal{O}(X)$, and the integers $n\geq 0$, $i\in \mathbb{Z}$. In the case, $S=\text{Spec}(k)$ for a subfield $k\hookrightarrow \mathbb{C}$, there are analogues of Betti and de Rham realization functors. There is also a version of Nori's Basic Lemma. This yields a tannakian category of exponential motives, which we denote by $M^{\text{exp}}(k)$. See [FJ20] for details of constructions. \\
This construction may extend to the relative case, leading to the tannakian category $M^{\text{exp}}(S)$ of relative exponential Nori motives.  
One can ask if the motivic theorem of the fixed part can be formulated in terms of the category $M^{\text{exp}}(S)$. A first obstruction is that there is no known theorem of the fixed part in the exponential category $M^{\text{exp}}(k)$ to start with.\\
 Because the (exponential) Gauss-Manin connections are no longer regular, one may replace the monodromy group by the \textit{differential Galois group} $G_{\text{diff},s}(M)$ associated to a relative exponential Nori motive $M/S$. Note that, in general,  the  characters of  $G_{\text{diff},s}(M)$ are no longer of finite order. One may ask whether, nevertheless, the natural sequence of algebraic groups is exact: $$ G_{\text{diff},s}(M)\rightarrow G_{\text{mot},s}(M) \rightarrow G_{\text{mot},0,s}(M)\rightarrow 0.$$  
 \\
  \vspace{20px}

\pagebreak

\begin{flushleft}
\Large\textbf{References}
\end{flushleft}
\phantomsection 
\addcontentsline{toc}{chapter}{Bibliography} 
\bibliographystyle{unsrtnat} 
\begin{flushleft}
[And92] Y. André, \textit{Mumford-Tate groups of mixed Hodge structures and the theorem of the
fixed part},  Compositio Math. \textbf{82} (1992), 1–24.
\end{flushleft}
\begin{flushleft}
[And96] Y. André, \textit{Pour une théorie inconditionnelle des motifs}, Publ. Math. IHÉS
\textbf{83} (1996), 5–49.
\end{flushleft}
\begin{flushleft}
[AB01] Y. André, F. Baladassarri, \textit{De Rham cohomology of differential modules on algebraic varieties}, Progress in Math.\textbf{189}, Brikhäuser, 2001.
\end{flushleft}
\begin{flushleft}
[And17a] Y. André, \textit{Groupes de Galois motiviques et périodes}, Séminaire Bourbaki.
Vol. 2015/2016. Exposés 1104–1119. Astérisque No. \textbf{390} (2017), Exp. No.
1104, 1–26.
\end{flushleft}
\begin{flushleft}
[And17b] Y. André, \textit{Normality criteria, and monodromy of
variations of mixed Hodge structures}, appendix C to B. Kahn, \textit{Albanese kernels and Griffiths groups}, arXiv: 1711.04335.
\end{flushleft}
\begin{flushleft}
[Ara13] D. Arapura, \textit{An abelian category of motivic sheaves}, Advances in Math. \textbf{233} 1 (2013), 135-195.
\end{flushleft}
\begin{flushleft}

[Ayo14] J. Ayoub. \textit{Periods and the conjectures of Grothendieck and Kontsevich-Zagier}, Eur. Math. Soc.
Newsl., (91):12–18, 2014.
\end{flushleft}

\begin{flushleft}
[BBD82] A. A. Beilinson, J. Bernstein, and P. Deligne. \textit{Faisceaux pervers}, In Analysis and topology
on singular spaces, I (Luminy, 1981), volume 100 of Astérisque, pages 5–171. Soc. Math.
France, Paris, 1982.
\end{flushleft}
\begin{flushleft}
[Bei87] A. Beilinson. \textit{On the derived category of perverse sheaves}, In K-theory,
arithmetic and geometry (Moscow, 1984–1986), volume 1289 of Lecture Notes
in Math., pages 27–41. Springer, Berlin, 1987.
\end{flushleft}

\begin{flushleft}

[CG14] U. Choudhury, M. Gallauer Alves de Souza. \textit{An isomorphism of motivic
Galois groups}, ArXiv e-prints, October 2014.
\end{flushleft}

\begin{flushleft}
[Del71] P. Deligne, \textit{Théorie de Hodge, II}, Publ. Math. IHÉS \textbf{40} (1971), 5–57.
\end{flushleft}

\begin{flushleft}
[DM82] P. Deligne, J. Milne,\textit{ Tannakian categories}, in Hodge Cycles, Motives, and
Shimura Varieties, Lect. Notes in Math. \textbf{900}, Springer, 1982, 101–228.

\begin{flushleft}
[DosSan15] J. P. Dos Santos, \textit{The homotopy exact sequence for the fundamental group
scheme and infinitesimal equivalence relations}, Algebr. Geom. \textbf{2} (2015), no.
5, 535–590.
\end{flushleft}
\begin{flushleft}
[EHS08] H. Esnault, P. Hai, X. Sun, \textit{On Nori’s fundamental group scheme}, Geometry and dynamics of groups and spaces, 377–398, Progr. Math. \textbf{265},
Birkhäuser, 2008.
\end{flushleft}

\begin{flushleft}
[HMS17] A. Huber, S. M\"{u}ller-Stach, \textit{Periods and Nori motives}, (with contributions of B. Friedrich and J. von
Wangenheim), Springer Ergebnisse \textbf{65}, 2017.
\end{flushleft}

\end{flushleft}
\begin{flushleft}
[Jos16] P. Jossen, \textit{On the relation between Galois groups and motivic Galois groups}, preprint, (2016). 
\end{flushleft}
\begin{flushleft}
[FJ18] J. Fres\'{a}n, P. Jossen, \textit{Algebraic cogroups and Nori motives}, arXiv: 1805.03906.
\end{flushleft}

\begin{flushleft}
[FJ20] J. Fres\'{a}n, P. Jossen, \textit{Exponential motives}, 2020, version July/20, manuscript available at http://javier.fresan.perso.math.cnrs.fr/expmot.fr
\end{flushleft}
\begin{flushleft}
[Kat70] N.M. Katz -- \textit{Nilpotent connections and the monodromy theorem: applications of a result of Turrittin,} IHES, \textbf{39,} (1970), pp. 175-232.
\end{flushleft}
\begin{flushleft}
[Kat72] N. Katz,\textit{ Algebraic solutions of differential equations (p-curvature and the Hodge filtration)}, Invent. Math.
\textbf{18} (1972), 1-118.
\end{flushleft}

\begin{flushleft}
[Kat82] N. Katz, \textit{A Conjecture in the Arithmetic Theory of Differential Equations}, Bull. Soc. math.
France \textbf{110}, 203–39.
\end{flushleft}

\begin{flushleft}
[Wat72] W. Waterhouse, \textit{Introduction to affine group schemes}, Graduate Texts in Mathematics, vol. \textbf{66} (Springer, Berlin, 1979).
\end{flushleft}
\vspace{10px}

\begin{flushleft}
LAGA, UMR 7539, CNRS, \textsc{University of Paris 13, University of Paris 8, Saint-Denis, France}\\
\textit{Email address}: mostaed@math.univ-paris13.fr,\;  amir.mostaed02@univ-paris8.fr
\end{flushleft}

\end{document}